# Efficient Mobility-on-Demand System with Ride-Sharing

Xianan Huang, Huei Peng

*Abstract*— An algorithm to cluster mobility-on-demand trips considering road network structure is developed in this paper. The benefits of our network partition algorithm are demonstrated in numerical simulations, showing that we can use fewer vehicles and can serve more customers with slightly longer wait time by including predicted future travel demand in trip assignment, compared with the benchmark reactive control policy.

## I. INTRODUCTION

Mobility-on-demand (MOD) services such as Uber and Lyft have brought significant changes, especially in urban areas with dense population. When multiple passengers share the same vehicle (e.g., Lyft Line and UberPOOL), the service can reduce the number of vehicles on the road and reduce congestion and energy consumption.

Today's MOD fleet management largely reacts to trip requests without utilizing predicted future supply and travel demand distribution. Continuous approximation [1] is used to study the dynamics of fleet and influence of large fleet to congestion. To control the fleet directly, algorithms such as mixed integer programming [2], heuristic [3] and graph based decomposition [4] methods demonstrated that current travel demand for taxis in New York city can be fulfilled with 15% of the existing fleet [5]. A privacy-preserving algorithm was developed [6] to protect the location information of passengers without incurring significant performance drop. However, the potential of the fleet is not fully utilized due to the nature of reactive control policy.

Knowledge of travel demand distribution plays a vital role in the control of MOD fleet. For carpool service with private cars, travel data can be used to identify optimal combined trips for carpooling and can reduce daily car mileage by 44% [7]. Intelligent transportation techniques such as connected automated vehicle provide richer information about travel demand and enable centralized coordination for the MOD fleet. Han et al. [8] showed that with driverless MOD fleet, the direct control approach is 29% more efficient compared with current price-based indirect control. For service provided by commercial fleet, travel demand distribution can be used to control the idling vehicles for rebalancing [8]–[10] to better meet future trip requests when carpooling is not allowed. A sampling-based algorithm is also proposed to control ride-sharing fleet using predicted future trip request information [11]. However, the travel location distribution is either characterized with clusters from geometric coordinate of locations [7] or grid-based discretization [9], neither of which takes the structure of the road network into consideration.

The research is supported by U.S. Department of Energy under the award DE-EE0007212.

The authors are with the Department of Mechanical Engineering, University of Michigan, Ann Arbor, Michigan 48109-2133, USA (e-mail: xnhuang@umich.edu, hpeng@umich.edu).

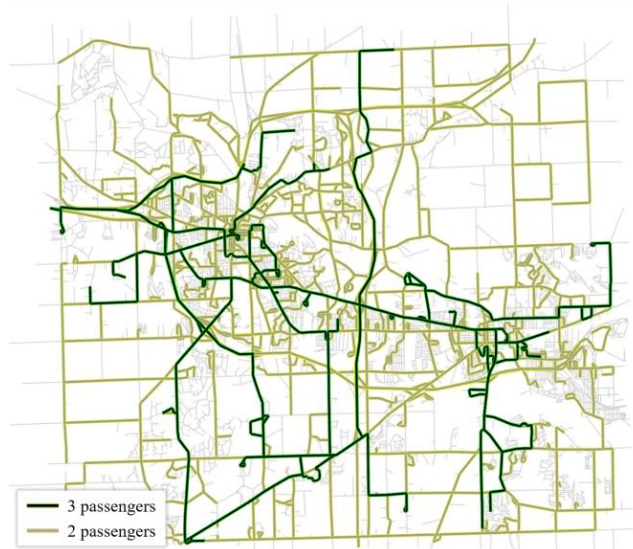

Fig. 1 Snapshot of shared trips of ride-share fleet operation with 900 vehicles serving 4% of travel demand during evening rush hour (17:00-18:00) in Ann Arbor

Since travel demand can be characterized as a random Origin-Destination variable on the road network, ignoring the underlying network structure can be problematic. To better describe the travel location distribution considering the structure of the road network, we propose an algorithm based on multidimensional scaling (MDS) [12] to project the locations on the road network onto an Euclidean space, and characterize the travel locations with Dirichlet Process Gaussian Mixture Model (DPGMM) [13]. The projection allows us to obtain better clustering results compared with the geometric coordinate-based methods. To utilize the demand distribution information for fleet management, we developed a fleet control algorithm based on the work in [5]. We propose a Kullback–Leibler (KL) divergence [14] regularization based control policy to balance current trip requests and future travel demand distribution. We assume that the demand distribution is known, and the fleet can be controlled directly to take assigned trips and to rebalance to be better prepared for future demands. It should be noted that we do not assume the future trips are known exactly, but their probability distribution is known. In our numerical study, travel demands are generated by POLARIS [15], a mesoscopic agent-based transportation model calibrated with data from the Safety Pilot Model Deployment (SPMD) project [16]. The calibration dataset consists of trip information from up to 2,800 vehicles since 2012.

The main contributions of this work are: 1) a travel location clustering algorithm considering the road network structure; 2) a ride-sharing fleet control policy considering future travel demand distribution for more efficient service.

The rest of this paper is organized as follows: Section 2 presents the proposed travel location clustering algorithm. Section 3 presents the formulation of demand regularized ride-sharing fleet optimization. Section 4 presents the simulation results. Conclusions and future work are given in Section 5.

## II. TRAVEL LOCATION CLUSTERING

To characterize travel location distribution with the road network structure taken into consideration, we model the original road network with a Euclidean space approximation. With this approximation, we characterize the distribution of travel locations with the Lebesgue measure.

### A. Multidimensional Scaling (MDS)

The MDS method is used to find the optimal Euclidian space that preserves pairwise distance in the network space. MDS can be formulated as an optimization problem defined as

$$min_{x_1,\ldots,x_N}\left(\frac{\sum_{i,j}(d_{ij}^p - \|x_i - x_j\|^p)^2}{\sum_{i,j}d_{ij}^{2p}}\right) \quad (1)$$

where $d_{ij}$ is the pairwise distance between points $i$ and $j$, $x_i, x_j \in \mathbb{R}^m$ are the vectors corresponding to point $i$ and $j$ in the projection space, and m is the dimension of the projection space. $p$ is the power transformation used by metric scaling, $N$ is the total number of projected points. Since the projection space is Euclidian, the approximated distance is

$$\|x_i - x_j\|^2 = \sum_{k=1}^{m}(x_{ik} - x_{jk})^2 \quad (2)$$

When $p$ is 1, the MDS is known as the classical MDS and can be solved with eigen-decomposition by transferring distance to inner product through double re-centering

$$G = -\frac{1}{2}\left(I_n - \frac{1}{n}\mathbf{1}\mathbf{1}^T\right)D\left(I_n - \frac{1}{n}\mathbf{1}\mathbf{1}^T\right) \quad (3)$$

where $D$ is the distance matrix, $D = \{d_{ij}\}$, $I_n$ is the identity matrix and $\mathbf{1}$ is the column vector with 1 as all its entries. With this transformation, vectors in the projection space can be obtained by eigen-decomposition of $G$, which gives

$$x_i^* = \sqrt{\lambda_i}u_i, i = 1 \ldots m \quad (4)$$

where $\lambda_i$ is the $i$-th largest eigenvalue of $G$, and $u_i$ is the corresponding eigenvector. When $p \geq 2$, the optimization problem can be solved using the steepest gradient method [12] where the solution of the classical MDS is used as the initial point for the numerical algorithm. In the following analysis, we use non-classical MDS with $p = 2$ to approximate the pairwise distance in the non-Euclidean road network space.

The distance matrix is obtained by calculating the pairwise lowest cost path distance between every pair of links in the capacity-normalized traffic network. The capacity-normalized traffic network is defined as a weighted directed graph with nodes associated with links of the original road network and edges associated with the movements. An edge from node $i$ to node $j$ exists if link $j$ is adjacent to link $i$ and if vehicle can travel from link $i$ to link $j$ (one-way road link is an example of when this is not the case). The weight of the edges is defined as

$$w_{ij} = \frac{1}{2}\left(\frac{l_i}{\bar{v}_i n_i} + \frac{l_j}{\bar{v}_j n_j}\right) \quad (5)$$

where $l_i, l_j$ are lengths of links $i$ and $j$, $\bar{v}_i, \bar{v}_j$ are the travel speeds, which can be the posted speed limits of the road links if no real-time traffic information is available, $n_i, n_j$ are lane numbers of the corresponding links. The graph is connected since there is no isolated links in the traffic network. The

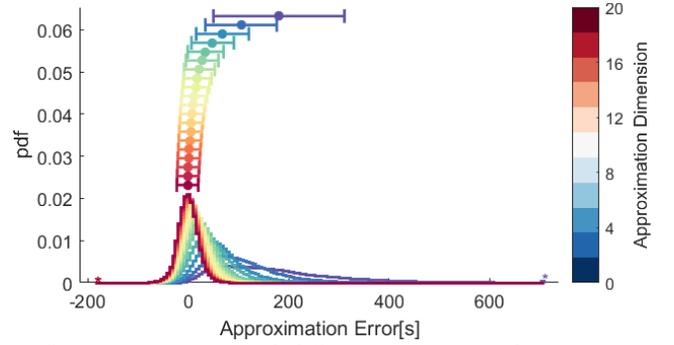
Fig. 2 Approximation Error with Different Approximation Dimension

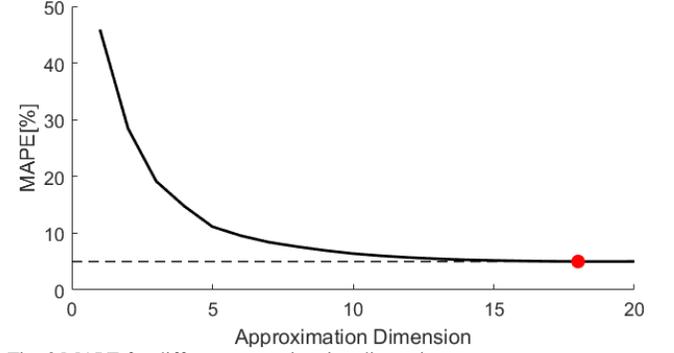
Fig. 3 MAPE for different approximation dimension

pairwise distance is solved using linear programming based on the Bellman inequality, which is the dual of Bellman-Ford algorithm [17] and can be solved efficiently with optimization solver such as Gurobi [18] which we used. The approximation performance with different projection space dimension is shown in Fig. 2, with mean and standard deviation marked with error bars.

In the following analysis, we choose 18 dimensions for the Euclidean space approximation, and the mean absolute percentage error (MAPE) converges as indicated by the red dot in Fig. 3. The MAPE for distance approximation is 5%, which is sufficient to preserve the pairwise distances of the original road network.

### B. Location Distribution Characterization

In the literature, travel location distribution is frequently characterized using Cartesian coordinate in the geometric space [7], which ignores road network structure information. With the approximation of link locations in the projected Euclidean space, we can analyze the location distribution in the projected space using the Lebesgue measure, which preserves the original distances in the network. In the following analysis, we assume that the origin and destination of each trip are sampled from the location distribution, and the union of origins and destinations is defined as locations of interests. Since we do not assume knowing the number of locations of interests, we use the DPGMM to model the random variable. DPGMM is a Bayesian nonparametric extension of the Finite Gaussian Mixture Model whose probability density function can be expressed by:

$$f_X(x) = \sum_{k=1}^{K}\pi_k f_{X,k}(x) \quad (6)$$

where $X$ is the random variable for travel locations, $f_X(x)$ is the overall density function, $\pi_k$ is the mixing coefficient for each component, $f_{X,k}(x)$ is the density for each component, which follows a multivariate Gaussian distribution. Instead of a fixed component number $K$, DPGMM assumes the model consists of infinite components, i.e. $K \to \infty$ in (6). With this method, not only the parameters for each mixture component but also the number of mixture components can be inferred from data. In this way, the locations of interests are modeled as the mean of each mixture component, and travel demand can be modeled as a multinomial distribution with the discretization achieved using the mixture model. The posterior of parameters of the DPGMM is inferred through collapsed Gibbs sampling, which is an approximate inference algorithm based on Markov Chain Monte Carlo (MCMC) sampling and known to be unbiased asymptotically compared with other approximate inference methods such as variational inference. The process is summarized as follows. Denote $c_i \in \{1, \ldots, K\}$ as the indicator variable of the component for each data point, which is a discrete random variable parameterized by $\pi$

$$P(c_i|\pi) \sim Discrete(\pi_1, \ldots, \pi_K) \quad (7)$$

For Dirichlet process model, $K \to \infty$. The parameters are modeled with their corresponding conjugate priors, i.e. Dirichlet distribution for $\pi$ and Gaussian-Wishart distribution for mean $\mu$ and covariance $\Sigma$ of each component.

$$P(\pi) = Dir(\pi|\alpha_0) = C(\alpha_0) \prod_{k=1}^{K} \pi_k^{\alpha_0 - 1} \quad (8)$$

$$P(\mu, \Sigma) = P(\mu|\Sigma)P(\Sigma)$$
$$= \prod_{k=1}^{K} N(\mu_k|m_0, \beta_0 \Sigma_k) W(\Sigma_k^{-1}|W_0, v_0) \quad (9)$$

where $\alpha_0, m_0, \beta_0, W_0, v_0$ are the hyperparameters. For simplicity, we denote $\{m_0, \beta_0, W_0, v_0\}$ the set of hyperparameters for the Gaussian-Wishart distribution as $\gamma$. The hidden variables include the indicator variable $c_i$ and the model parameters $\pi, \mu, \Sigma$. At each step of collapsed Gibbs sampling, we sample $c_i$ conditional on the rest of the data points and random variables from the posterior

$$p(c_i = k|c_{-i}, x, \alpha_0, \gamma)$$
$$\propto p(c_i = k|c_{-i}, \alpha_0) p(x|c_{-i}, c_i = k, \gamma) \quad (10)$$

where $c_{-i}$ is the set of indicator variables for other samples except $i$, $c_{-i} = \{c_j, j \neq i, j \in \mathbb{N}, 1 \leq j \leq N\}$, $N$ is the sample size for the entire dataset. Since the prior for other parameters are well-defined, the inference can be carried out in a closed form. Thus, no sampling is required to obtain the posterior of $\mu$ and $\Sigma$ once $c_i$'s are sampled for all data points.

The likelihood term can be obtained in a closed form from the Gaussian-Wishart distribution, and the prior term can be defined by the Chinese Restaurant Process (CRP). The resultant cluster assignment follows the pattern that the probability of a new sample belonging to a cluster is proportional to the number of samples already in the cluster.

$$p(c_i = k|c_{-i}, \alpha_0) = \begin{cases} \dfrac{N_{-i,k}}{N + \alpha_0 - 1} & \text{If } k \leq K \\ \dfrac{\alpha_0}{N + \alpha_0 - 1} & \text{If } k = K + 1 \end{cases} \quad (11)$$

where $N_{-i,k}$ is the sample size of data belong to cluster $k$ for other samples except $i$, $K$ is the current number of clusters already realized. In this way, as the sample size $N$ goes to infinity, the number of clusters can go to infinity, indicating that the model is more complex with more samples acquired.

DPGMM is used to identify the clusters for travel locations in the projected Euclidean space. The clustering result from DPGMM is used to partition the network and define the sample space for travel demand distribution, which is defined on the origin-destination region pairs.

### III. TRAVEL DEMAND REGULARIZED ASSIGNMENT

Our fleet control algorithm is based on the graph decomposition method proposed in [5]. The algorithm can solve the trip matching and routing problem for ride-sharing for thousands of vehicles and customers fast enough for real-world implementation. We further improve the algorithm to take knowledge of future travel demand distribution into consideration.

#### A. Real-Time Ride-share Trip Assignment

As a start point, we reproduce the work in [5] by assuming the road network is static and solving all optimal routes considering only travel time offline. Including dynamic road network information is done later. The trip assignment algorithm is based on a shareability graph. The graph is defined as undirected graph with nodes defined as customers and vehicles. An edge exists between two customers if a vehicle can depart from the origin of one of the customers and fulfill the travel demands of both customers without violating travel time constraints. An edge exists between a vehicle and a customer if the demand can be served by the vehicle without violating travel time constraints. Then a necessary condition for a trip to be feasible is that the customers of the trip can form a clique with one vehicle present in the shareability network. A clique is a subgraph such that every node is connected to every other node within the same clique. It's noted that the cliques do not need to be maximum cliques in the shareability graph. The cliques in a graph can be found with Bron-Kerbosch algorithm [19] with worst case time complexity $O(dn3^{d/3})$ where $n$ is the number of nodes and $d$ is degeneracy of the graph, which is a measure of sparseness. In this way, instead of evaluating cost of trips for every possible combination of customers and vehicles, one can solve single-vehicle-multiple-customer problems for every clique.

Trip scheduling for each clique is a traveling salesman problem with pickup and delivery. The problem can be solved with multiple algorithms. If the number of customers is small, (e.g., less than 5), the exact solution can be found by Dynamic Programming in less than 1 sec on a standard desktop computer. Heuristic based algorithms such as T-share [20] can be used to find the solution if the problem size is large.

After all feasible trips were found through solving the scheduling problem for all cliques, the optimal trip assignment problem can be formulated and solved through Integer Linear Programming (ILP). In this Section, we briefly summarize the formulation from [5] and the additional regularization term for demand distribution is presented in next section.

The cost for each customer consists of wait time and delay time. Wait time is defined as time between the customer travel request and time of pickup. Delay time is defined as the difference between planned travel time and the shortest travel time after pickup, which is from the fastest path solution from origin to destination. The cost of a trip is defined as wait time plus delay time for all customers, denoted as $c_t^i$ for trip $i$. The states of the system are $\delta_t$ which is the indicator variable for trip/clique and $\delta_c$ which is the indicator variable for a customer. If at an assignment instant, there are $m$ feasible trips from TSP step and $n$ customers, then $\delta_t = \{\delta_t^i \in \{0,1\}, i \in \mathbb{N}, 1 \le i \le m\}$ and $\delta_c = \{\delta_c^i \in \{0,1\}, i \in \mathbb{N}, 1 \le i \le n\}$. $\delta_t^i$ is 1 if trip $i$ is selected and $\delta_c^i$ is 1 if customer $i$ is assigned. The objective function is

$$\sum_{i=1}^{m} c_t^i \delta_t^i + \sum_{i=1}^{n} D(1 - \delta_c^i) \quad (12)$$

where $D$ is the penalty for unserved customers. The constraint for vehicle is that each vehicle can only serve one trip

$$\sum_{i=1}^{m} a_j^i \delta_t^i \le 1, \forall j \quad (13)$$

where $a_j^i$ is the indicator variable for vehicle $j$ and trip $i$, $a_j^i = 1$ if vehicle $j$ can serve trip $i$. The constraint for customer is that a customer is either assigned or ignored

$$\sum_{i=1}^{m} b_j^i \delta_t^i + (1 - \delta_c^j) = 1, \forall j \quad (14)$$

where $b_j^i$ is the indicator variable for customer $j$ and trip $i$, $b_j^i = 1$ if customer $j$ can be served by trip $i$. With linear constraints and the objective function, the trip assignment problem is an integer linear programming. For online optimization, we follow [5] to keep a pool of customers and a customer is removed from the pool if it's picked up by vehicle or cannot be served within the time constraint. If a customer is ignored, a vehicle from the idling fleet is assigned to serve the vehicle with minimum wait time as the objective.

### B. Travel Demand Regularization

The algorithm developed in [5] is limited by its reactive nature. The follow-up work [11] uses a sampling method to keep track of demand distribution. However, with 400 samples of future demands the computation time increased significantly. Since the dimension of the demand distribution is $K^2$ where $K$ is the number of partitions of the network, the sample size needed to characterize the distribution can be large. To balancing the fleet assignment and potential travel demand, we propose a regularization based on KL divergence between vehicle empty space distribution and travel demand distribution. We model the vehicle space distribution and travel demand distribution as a multinomial distribution defined on origin-destination region pairs, which forms a finite discrete sampling space. For an assignment instant, the distribution of vehicle space can be obtained from maximum likelihood estimation and the demand distribution is assumed to be known. The KL divergence is an asymmetric measurement of difference between two distributions, and is defined as

$$KL(q \parallel p) = \sum_i q_i \log \frac{q_i}{p_i} \quad (15)$$

where $q_i$ is probability for OD pair $i$ in the demand distribution and $p_i$ is probability for OD pair $i$ in the vehicle space distribution. The vehicle space distribution can be estimated by the trip plan from the TSP step

$$\hat{p}_i = \frac{\sum_j s_j \delta_{v,j}^i}{\sum_i \sum_j s_j \delta_{v,j}^i} \quad (20)$$

where $s_j$ is the number of available space for vehicle $j$, $\delta_{v,j}^i$ is the indicator variable for vehicle $j$, the value is 1 if the vehicle is traveling between OD pair $i$. The denominator is the total available space for the fleet, which can be approximated by the total space of fleet minus the number of customers, which is a constant. Thus, we can focus on the numerator in the optimization. The demand distribution is assumed to be known and invariant to trip assignment

$$KL(q \parallel p) = -\sum_i q_i \log p_i + const$$
$$= -\sum_i q_i \log \sum_j s_j \delta_{v,j}^i + const \quad (16)$$

Since the log function is concave, we can apply Jensen's inequality [21] to the nonlinear term to get the upper bound of the KL divergence

$$KL(q \parallel p) \le -\sum_i q_i \sum_j \delta_{v,j}^i \log s_j + const \quad (17)$$

This approximation transforms the nonlinear KL divergence to a linear function for indicator variable $\delta_{v,j}^i$. In this way, we can put the regularization term into the objective function as an additional linear term. With the additional indicator variable for the vehicles, the constraint for vehicle assignment changes to an equality constraint

$$\sum_{i=1}^{m} a_v^i \delta_t^i + \sum_j (1 - \delta_{v,j}^i) = 1, \forall j \quad (18)$$

And the new objective function is

$$\sum_{i=1}^{m} c_t^i \delta_t^i + \sum_{i=1}^{n} D(1 - \delta_c^i)$$
$$- w_R \sum_i q_i \sum_j \delta_{v,j}^i \log s_j \quad (19)$$

where $w_R$ is the weighting parameter for regularization. With the approximation applied, the objective function and constraints are linear, thus the problem is still an ILP which can be solved efficiently. It should be noted that generally integer programming is hard to solve due to the combinatorial nature. However, to solve an integer programming efficiently is beyond the scope of this research and we use Gurobi to solve the optimization problems.

## IV. RESULTS AND DISCUSSION

In this section, we present the traffic network partition results using the proposed algorithm. The traffic demand is generated using POLARIS. In this study, we focus on demand generated during the evening rush hour (17:00-18:00) on weekdays. However, the algorithms developed can be extended to deal with time-varying demand distribution which can be modeled as a piece-wise constant function. We first present the results of road network partition, then using a numerical simulation to demonstrate the demand distribution regularized fleet control.

### A. Road Network Partition

To evaluate the performance of clustering algorithm, we used the origins and destinations of trips generated from 17:00 to 18:00 on a weekday by POLARIS. The heat map of the origins and destinations is shown in Fig. 4(a). The network partition result from our proposed algorithm using DPGMM in the projected Euclidean space is shown in Fig. 4(b). 30

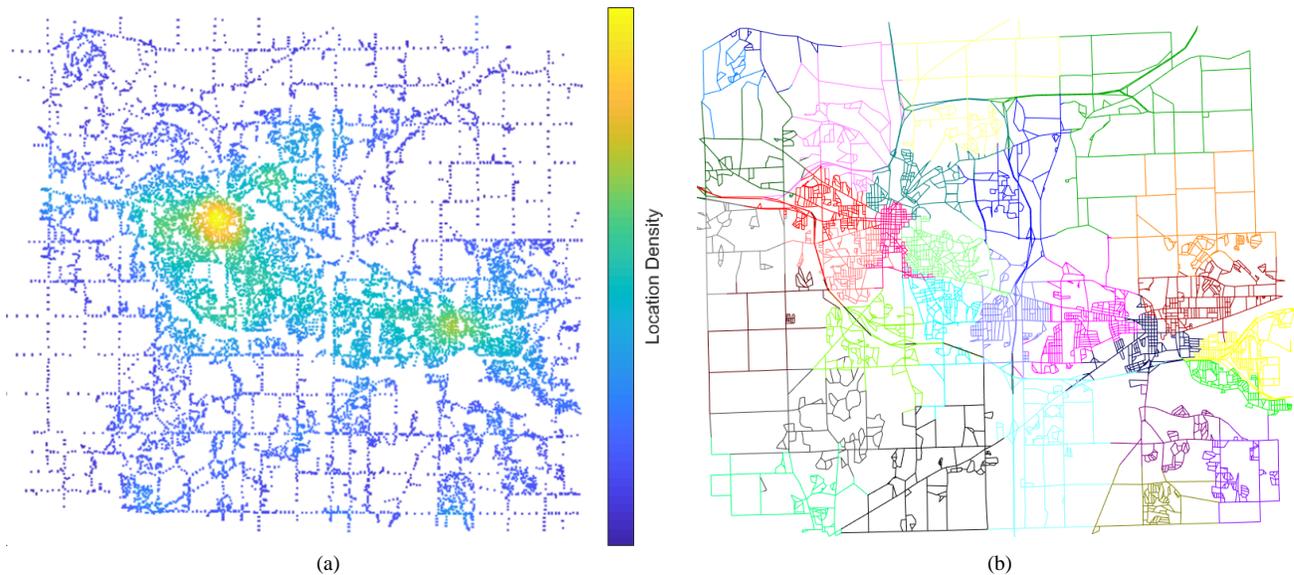

Fig. 4 (a) Locations of interest heatmap generated by POLARIS simulation from 17:00 to 18:00 on weekdays in Ann Arbor (b) network partition using the proposed algorithm with DPGMM in the projection space for POLARIS generated demand components were identified with the Bayesian nonparametric algorithm, with different clusters in the original network space indicated with different color. With the partition result, travel demand can be described with a discrete random variable following multinomial distribution with sampling space size 900.

The benchmarks for traffic network partition are clustering algorithms applied in the geometric Cartesian coordinate space including density based algorithm GMM and distance based algorithm k-means [14]. Since the traffic network partition we are interested in is clustering data in the network space instead of studying the connectivity of the network itself, the community detection algorithms [22] are out of scope for our evaluation. Since the component number needs to be specified for k-means, to make the evaluation a fair comparison, instead of using Bayesian nonparametric algorithm to identify component number for GMM and our proposed algorithm, we use the Expectation-Maximization (EM) algorithm to identify parameters for the mixture models. The performance metric we selected is mean travel time to cluster center, with travel time for each road section generated from POLARIS during the studied hour. Since our objective for clustering is to approximate the infinite dimensional travel demand distribution with finite dimensional discrete distribution, short within-cluster travel time is desirable so the locations within each cluster can be approximated with the cluster center. The average travel time to cluster center for different cluster algorithms are shown in Fig. 5, with right axis showing average travel time reduction with our algorithm compared with k-means in blue. The histogram for mean travel time to cluster center for different clusters with 30 mixture components is shown in Fig. 6. Since GMM in the Cartesian coordinate failed to generate enough clusters, the histogram is not included.

As shown in Fig. 5, GMM with geometric coordinate cannot fit model with more than 13 components from the data and more components would result in ill-conditioned covariance matrices during iterations of EM algorithm. However, this doesn't mean that more clusters would overfit the data since the actual cluster number can be estimated with our Bayesian nonparametric algorithm, which is 30. The failure to identify additional mixture components is due to the lack of network structure information. The results indicate that with the right number of mixture component, our algorithm can reduce in-cluster travel time by more than 10% compared with clusters in geometric coordinate as shown in Fig. 5. Also, as indicated in Fig. 6, with the same number of clusters, our algorithm have shorter travel time for all clusters compared with the benchmark, and the variance for distribution is smaller, with standard deviation 11.82 sec compared with 16.88 sec indicating that the network is partitioned more uniformly compared with the benchmark. In the next section, we use a numerical simulation to demonstrate that a good partition can benefit fleet control.

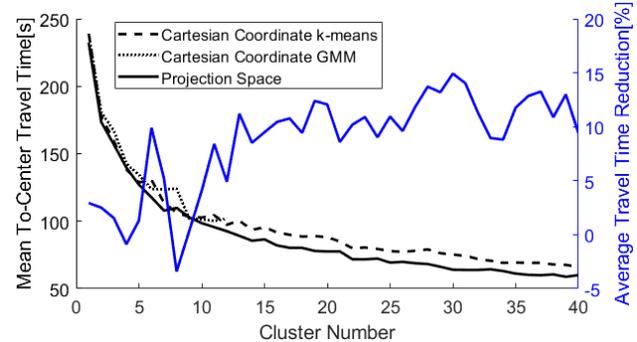

Fig. 5 Mean travel time to cluster center for all clusters using different cluster algorithm

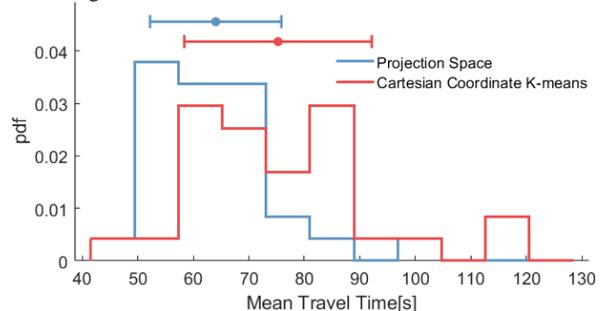

Fig. 6 Mean travel time for each cluster with cluster number=30

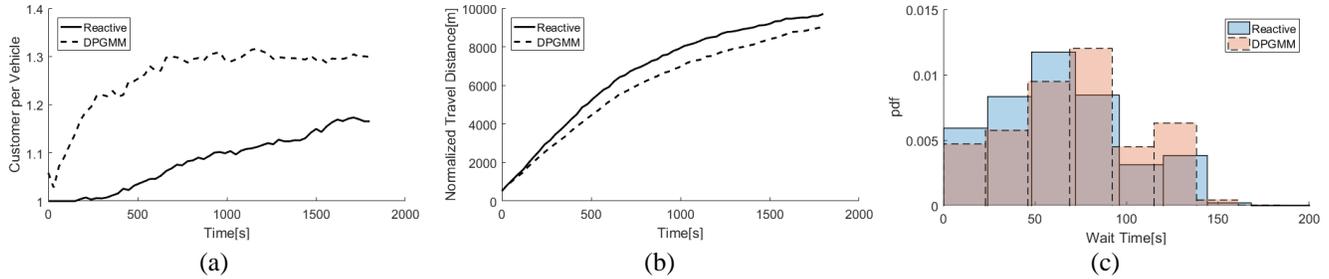

Fig. 7 Performance Comparison for Benchmarks: (a) Average planned customer number per vehicle; (b) Total travel distance normalized by number of served customers; (c) Customer wait time histogram for on time planned customer

## B. Numerical Results for Ride-Sharing Control

We randomly selected 4% of the trips generated during the studied time as demand for the shared mobility fleet. The trip generation rate is 35~40 new trips every 30 sec and we follow the re-optimization strategy every 30 sec from [5]. The simulation period for our study is 30 min. We fix the fleet size at 900 and the vehicle capacity is 4. The wait time constraint is 2 min, and the delay time constraint is 4 min. The benchmark algorithm is the reactive policy. Our proposed algorithm is denoted as DPGMM. We assume the fleet is perfectly balanced initially, achieved through sampling the vehicle initial locations randomly from the trip origins. The pairwise travel times between all pairs of roads are calculated offline using the average speed from POLARIS during the studied period. The performance metrics are shown in Fig. 7.

During the simulated period, all customers are served. As shown in Fig. 7(a), our proposed algorithm results in more customers per operating vehicle, indicating more efficient operation of the fleet. Besides that, as shown in Fig. 7(b), with normalized distance defined by fleet total travel distance normalized with number of served customers, the regularization also decreased the normalized travel distance of the fleet, which reduces the operation cost for the service provider. With the additional regularization term, a larger idling fleet is kept during the operation before all vehicles are assigned, meanwhile same number of customers are served with the proposed algorithm, indicating the fleet is utilized more efficiently. Since the regularization is achieved through weighted sum of KL divergence cost and the travel time cost, the regularization term results in an increase for mean customer wait time by 12%. However, it should be noted that both algorithms can serve 87% of customers within the time constraints. The results also indicate that the commonly used static network assumption is too strong. If the static network assumption holds, all trips should satisfy the wait time constraint, which is not the case as shown in Fig. 7 (c), indicating that with 900 vehicles in operation, the travel time changes in the network.

## V. CONCLUSION AND FUTURE WORK

We proposed a road network partition algorithm based on MDS and DPGMM. The partition is used to discretize the travel demand distribution. It can be integrated with fleet control strategy to utilize the fleet more efficiently with our proposed regularized assignment algorithm. However, currently we still have the static traffic network assumption, which results in performance degradation in our dynamic simulation. For our next step, we will improve the robustness of ride-sharing control algorithm to include traffic dynamics in the control strategy design.